\documentclass[11pt,letterpaper]{article}

\usepackage[T1]{fontenc}
\usepackage{lmodern}
\usepackage[margin=1in]{geometry}
\usepackage{amsmath,amssymb,booktabs,enumitem,microtype}
\usepackage[hidelinks]{hyperref}
\usepackage{xurl}

\hypersetup{
  pdftitle={A vertex operator algebra realization of the Haagerup double},
  pdfauthor={Andrew Riesen},
  pdfsubject={Research announcement}
}
\setlist[enumerate]{leftmargin=*,itemsep=4pt,topsep=4pt,parsep=2pt}
\setlist[itemize]{leftmargin=*,itemsep=1pt,topsep=2pt}
\newtheorem{theorem}{Theorem}
\newcommand{\A}{\mathcal A}
\newcommand{\g}{\mathfrak g}
\newcommand{\HH}{\mathcal H}
\newcommand{\CC}{\mathcal C}
\newcommand{\DD}{\mathcal D}
\newcommand{\one}{\mathbf1}
\DeclareMathOperator{\Com}{Com}
\DeclareMathOperator{\Rep}{Rep}
\DeclareMathOperator{\Hom}{Hom}

\title{The Haagerup VOA at $c=8$}
\author{Andrew Riesen}
\date{October 9th 2026}
\begin{document}
\maketitle

\begin{abstract}
We announce an explicit construction of a unitary, strongly rational
vertex operator algebra $U$ of central charge $8$ inside
$V_{E_6\oplus A_2}$, with a ribbon equivalence
$\Rep(U)\simeq Z(\mathcal H_3)$ realizing the Haagerup double.
The vertex algebra is generated by six affine currents and three
specified weight-two fields; its weight-one Lie algebra is
$\mathfrak{sl}_2\oplus\mathfrak{sl}_2$, at levels $1$ and $39$.
We outline the argument, which combines exact lattice calculations,
an all-weight strong-generation induction and $C_2$-cofiniteness,
semisimple descent for a full corner of Zhu's algebra, modular
Galois symmetry, and holomorphic $E_8$ extensions.

\end{abstract}

The point of this article is to announce that, using AI, the author
has found a strongly rational VOA of central charge $c=8$ which
realizes the double of the Haagerup fusion category. More precisely,
the following result has been shown through extensive calculations
by AI:

\begin{theorem}
There exists a unitary, strongly rational vertex operator algebra $U$
of central charge $8$, realized as a conformal subalgebra
\[
U\subset V_{E_6\oplus A_2}\subset V_{E_8},
\]
together with a $\mathbb C$-linear ribbon equivalence
\[
\Rep(U)\simeq Z(\mathcal H_3),
\]
where $\mathcal H_3$ is the Haagerup fusion category. The weight-one Lie
algebra of $U$ is $\mathfrak{sl}_2\oplus\mathfrak{sl}_2$, with affine
levels $1$ and $39$.
\end{theorem}

In this writeup, we give an outline of the proof, the techniques used,
and the required calculations. Due to the unusual manner in which
this was discovered, we believe it is best to first disclose how
the discovery happened and the role the author has played.
\section*{AI Disclosure}

The text of this announcement was written by the author, with AI used for editing. The working manuscript and certificates combine AI-generated material with the author's own contributions. The author takes full responsibility for this announcement, the working manuscript, and the accompanying certificates.
\section*{AI Methodology}

The main AI tools used have been GPT Pro 6, Astra 6, and Sol 6.1.
Here is the order of events leading to the result. For the sake of
brevity, we omit the many prompts required to reorient the AI
properly.

\begin{enumerate}
    \item The author asked Astra 6 to use all available methods
    to find candidate VOAs which realize $Z(\mathcal H_3)$.
    After 3 hours and 22 minutes, Astra 6 provided a promising
    explicit conformal subVOA
    \[
    U\subset V_{E_6\oplus A_2}
    \]
    that satisfied several nontrivial tests:
    \begin{enumerate}
        \item $\dim U_2=120$ and $\dim U_3=956$;
        \item $U$ satisfied several compatibility conditions
        required by vector-valued Jacobi forms;
        \item $U$ had nilpotence relations in the $C_2$-quotient.
    \end{enumerate}

    \item The author prompted GPT Pro 6 to prove that $U$
    is finitely strongly generated and $C_2$-cofinite.
    Pro 6 claimed to have proved both.

    \item The author prompted Astra 6 to show that $U$ has
    $12$ simple modules and to determine their fusion rules.
    The AI claimed that it had rigorously constructed six
    non-isomorphic $U$-modules using lattice constructions.
    In a subsequent conversation, GPT Pro 6 found a subVOA arising from
    a simple current extension $\mathcal A^+\subset U$
    and a full idempotent $e_Z$ reducing the problem to
    an easier algebra called the corner algebra.

    \item The author prompted Astra 6 to continue in another
    new conversation. The AI claimed that, with the assistance
    of code, it had continued the previous work and shown
    that $U$ is strongly rational, has $12$ simple modules,
    and has a Zhu algebra of dimension $36{,}625$.

    \item The author then discussed with Pro 6 potential proof
    strategies for showing the categorical equivalence.
    The AI pointed out that $U\subset V_{E_8}$, and the author
    realized that the corresponding extension would give a
    Lagrangian algebra in $\Rep(U)$.

    \item The author prompted Sol 6.1 to establish a braided
    equivalence between $\Rep(U)$ and $Z(\mathcal H_3)$.
    The AI thought for 1 hour and 6 minutes and provided a proof.

    \item The author prompted Sol 6.1 to upgrade the braided
    equivalence to a ribbon equivalence. The AI did so.
\end{enumerate}

After this point, the author prompted the AI to create a zip file
consisting of a manuscript, the data and code used in the
calculations, and verifiers. The AI provided this material.

Through sustained engagement with the arguments and their
presentation, the author developed an understanding of the proof
strategy and its principal steps. This work increased the author's
confidence in the proposed proof to the point at which they judged
that a public announcement was warranted. A central contribution of the author has been to refine and hone
the arguments, notation, and exposition, preparing the manuscript
and computational material for human understanding, verification,
and further mathematical contributions. The author also directed
the search and guided the proof strategy, including recognizing the
role of Lagrangian algebras in the categorical identification.

This work remains ongoing, and the author would greatly appreciate
community feedback on the arguments and calculations, their
presentation, and appropriate attribution.

\section{Structure of the Project and Calculations}

The detailed calculations are given in the accompanying working manuscript and certificate package. The manuscript is being refined for a forthcoming second version. The GitHub repository for the project is:
\begin{center}\url{https://github.com/Doctorfm98/Haagerup_VOA}
\end{center}
The project is being organized into two linked components:
\begin{enumerate}
\item \textbf{The manuscript,} with appendices explaining the calculations and their inputs and outputs, indexed by labels such as \textsf{C01} and \textsf{C20}.
\item \textbf{Computations and certificates,} organized by the same labels. Each guide identifies its inputs and witnesses, implementation, verification entry point, and instructions for reproduction and human inspection.
\end{enumerate}

The reference package is \path{Haagerup_VOA_Post.zip}. Its working manuscript is \path{Haagerup_VOA_Post.pdf}, with editable source starting at \path{manuscript/main.tex}. Certificate labels \textsf{C01--C21} refer to the twenty-one entries in \path{CERTIFICATES.md}, which link to the corresponding guides under \path{certificates/}. The manuscript's calculation registry in Subsection~R.2 uses the same labels. All package paths are relative to the extracted archive root; unqualified section, proposition, and appendix numbers refer to this working manuscript.
\section*{The explicit candidate}
Put
\[
W=V_{E_6\oplus A_2},\qquad
\A=L_1(\mathfrak{sl}_2)\otimes L_{39}(\mathfrak{sl}_2).
\]
The candidate is the vertex subalgebra
\[
U=\langle\A,e,f,v\rangle_{\mathrm{VA}}\subset W.
\]
The first affine factor comes from a root \(A_1\) in \(E_6\). Orthogonal to it is a chosen \(A_2\oplus A_5\) root system; the diagonal of its principal \(\mathfrak{sl}_2\) triples has level \(4+35=39\). The six currents and \(e,f,v\) are the nine vertex generators; strong generation is a separate assertion.

Cartan-refined character constraints and the ambient affine decomposition guided the choice of two coset fields and a weight-two highest-weight space of type \((1,13)\). The formulas use the root embedding, cocycle, and ordered bases fixed in Section~2 and Appendices~B--E of the working manuscript.

Write
\[
C_W=\Com(\A,W),\qquad \omega_C=\omega_W-\omega_{\A}.
\]
In the specified basis \(b_0,\ldots,b_3\) of the four-dimensional space \((C_W)_2\), take
\[
\begin{aligned}
\omega_C&=\frac{8b_0+9b_1+2b_2+2b_3}{41},\\
e&=\frac{2b_0+6b_1}{31}+\frac{2-r}{93}b_2+\frac{2-2r}{93}b_3,
\qquad f=\omega_C-e,\qquad r=\sqrt{13}>0.
\end{aligned}
\]
These are commuting Virasoro vectors, of central charges \(70/31\) and \(2400/1271\). Let \(u_i=e^{\beta_i}\) be the specified exponential basis of the three-dimensional weight-two affine highest-weight space of type \((1,13)\), and choose
\[
v=\frac{-r+3q}{20}u_1+\frac{r+3q}{20}u_2+u_3,
\qquad q=\sqrt{-3}=i\sqrt3.
\]
The coset eigenvalues \(e_{(1)}v=13v/31\) and \(f_{(1)}v=180v/1271\) select a two-dimensional space, so the displayed coefficients specify an additional choice of line. The data define a \(K\)-form of \(U\), where
\[
K=\mathbb Q(r,q)=\mathbb Q(\theta),\qquad
\theta=r+q,\qquad \theta^4-20\theta^2+256=0.
\]
Since \(\omega_{\A}+e+f=\omega_W\in U\), the inclusion is conformal and \(c(U)=8\). The definition and these finite checks are recorded under \textsf{C01}.

\section*{Strong rationality and the Zhu algebra}
The first objective is to prove that \(U\) is simple, self-contragredient, of CFT type, rational, and \(C_2\)-cofinite. These conditions make \(\Rep(U)\) a modular tensor category. The proposed route is as follows.

\subsection*{1. Simplicity, unitarity, and self-contragredience}
A transported lattice conjugation preserves \(U\), with an explicit identity placing the conjugate of \(v\) back in \(U\) (\textsf{C02}; Section~3 and Appendix~G). Its positive invariant Hermitian form restricts to \(U\); the associated nondegenerate invariant bilinear form gives self-contragredience, and a least-weight ideal argument gives simplicity. The inherited grading gives CFT type.

\subsection*{2. Strong generation and \(C_2\)-cofiniteness}
Low-weight mode calculations and Jacobi relations supply a proposed strong generating set in weights at most four (\textsf{C03--C06}; Sections~4--5 and Appendix~H). A triangular mode induction extends the finite closure checks to every weight. This is essential: low-weight ranks alone do not prove strong generation.

Strong generation makes \(R_U=U/C_2(U)\) a finitely generated commutative Poisson algebra. Its nilradical is a Poisson ideal in characteristic zero. Affine integrability, current covariance, and exact square identities for \(e,f\) put the weight-one and weight-two generators in this ideal (\textsf{C07}; Section~5 and Subsection~H.9). The higher generators are built from them by translations, products, and zeroth products, so are also nilpotent. Finitely many nilpotent algebra generators then make \(R_U\) finite-dimensional.

\subsection*{3. Semisimplicity through a full corner and a path algebra}
The aim is to prove semisimplicity of the whole Zhu algebra before exhausting the modules.

\begin{enumerate}[label=(\alph*)]
\item \textbf{An affine extension.} An explicit lattice product places the simple current \(J_{(1,39)}\) inside \(U\), giving
\[
\A^+=\A\oplus J_{(1,39)}\subset U.
\]
Its twenty ordinary simples have known lowest spaces. Each is irreducible under \(\g=\mathfrak{sl}_2\oplus\mathfrak{sl}_2\), except for
\[
T_{20}=V(0,20)\oplus V(1,19),\qquad \dim T_{20}=21+40=61,
\]
using the affine section's even labels (\textsf{C10}; Section~7 and Appendix~J). Equality of inherited and affine gradings and translations gives the unital map of Proposition~7.2:
\[
\Phi:A(\A^+)\longrightarrow A(U).
\]
For the matrix and quiver coordinates, the even affine block labels \(j\) are renumbered as \(t=j/2\in\{0,\ldots,19\}\), so the exceptional block above is vertex \(10\). In these coordinates, \(T_t\) denotes the lowest space of block \(t\).

\item \textbf{The full corner.} For the images \(E_{ij}^{(t)}\) of the affine block matrix units, put
\[
e_t=E_{11}^{(t)},\qquad e_{\rm Z}=\sum_t e_t,
\qquad D=e_{\rm Z}A(U)e_{\rm Z}.
\]
The identity
\[
1=\sum_{t,i}E_{i1}^{(t)}e_{\rm Z}E_{1i}^{(t)}
\]
proves fullness, including when some blocks vanish under \(\Phi\). Thus \(D\) is Morita equivalent to \(A(U)\).

\item \textbf{Scalar matrices and coefficient arrows.} The same matrix units give an algebra realization
\[
A(U)\cong\bigoplus_{t,v}
\Hom_{\mathbb C}(T_t,T_v)\otimes e_vDe_t,
\]
with multiplication by matrix composition and corner-coefficient multiplication. The transported current action is on the known matrix factor and trivial on the coefficient factor. This does not require the embedded corner to be stable under current commutators.

The all-weight reduction gives twenty-two Zhu generators: six currents and sixteen highest-weight classes \([z_0],\ldots,[z_{15}]\) (Proposition~6.2 and Appendix~A). For each admissible irreducible constituent, a fixed normalized highest-weight matrix extracts a unique coefficient
\[
\xi_{t\to v}^{\lambda,\mu}(z_a)\in e_vDe_t.
\]
These coefficients and the vertex idempotents generate \(D\), defining a quiver over \(K\). Section~8 and Appendix~L fix the numerical matrices and coefficient extraction; Definition~9.1 specifies the quiver and its path algebra. Paths are composable arrow products.

\item \textbf{From lattice identities to path relations.} Three current-equivariant lattice-state maps compare the following expressions in the computed spaces \(S_2\) and \(S_3\):
\begin{center}
\begin{tabular}{@{}lr@{}}
\toprule
Relation family & Highest-weight families\\
\midrule
Symmetric normal products of weight-two states & 66\\
Translations, current products, and alternating zeroth products & 408\\
Mixed weight-four relations involving weight-three states & 446\\
\bottomrule
\end{tabular}
\end{center}
For each map, a kernel vector is an expression that vanishes as a lattice state. The three kernels give \(66+408+446=920\) highest-weight VOA relation families (\textsf{C11}; Section~9 and Appendix~K); applying the current lowering operators gives the remaining relations in each family. To obtain a relation in Zhu's algebra, we project the state identity to \(A(U)\), replacing translations by their weight-dependent Zhu classes and retaining the lower-mode correction terms. We then express each generator through the fixed numerical matrices and its corner coefficients, and compare matrix entries. The resulting equations are polynomials in composable arrows, so they define path-algebra relations. A single highest-weight VOA relation can yield several such polynomials; \(920\) counts the original relation families, not the final path relations.

Together with conformal centrality, these give an upper presentation. Exact consequences and eliminations give a smaller presentation \(S_K\) with
\[
S_K\otimes_K\mathbb C\twoheadrightarrow D.
\]
This surjection suffices for semisimplicity; identifying the complete relation kernel is unnecessary at this stage.

\item \textbf{The integral model and its special fiber.} Set
\[
R=\mathbb Z[\theta]_{(10009,\theta-2914)},\qquad
\mathfrak m=(10009,\theta-2914)R.
\]
This is a DVR with fraction field \(K\) and residue field \(\mathbb F_{10009}\). The exact eliminated rules have coefficients in \(R\) and define \(S_R\) with
\[
S_R\otimes_R K=S_K.
\]
Integral shortening of every three-arrow path shows that \(1{,}081\) specified paths of length at most two span \(S_R\). This proves finiteness over \(R\) before specialization. The special fiber \(S_R/\mathfrak mS_R\) has dimension \(546\), a nondegenerate regular trace pairing, and an explicit split matrix decomposition with block sizes
\[
1,1,1,1,2,3,4,6,9,10,10,14.
\]
Vertex exclusions and exact quadratic reduction are recorded under \textsf{C12--C13} (Subsection~K.4), integral finiteness under \textsf{C14} (Subsections~K.5.1--K.5.2), and the special-fiber structure under \textsf{C15} (Subsections~K.5.3--K.5.4). The good-prime arithmetic is recorded under \textsf{C17} (Subsection~R.1).

\item \textbf{Descent and rationality.} Removing \(R\)-torsion and lifting surviving blocks over the completed DVR proves that \(S_K\) is absolutely semisimple, with \emph{at most} twelve absolute simple blocks. Blocks may disappear through torsion or through the quotient onto \(D\); their survival is proved below.

Consequently \(D\) is semisimple, and fullness transfers semisimplicity to \(A(U)\). With simplicity, self-contragredience, CFT type, and \(C_2\)-cofiniteness already established, McRae's rationality criterion~\cite[Corollary~5.19]{McRae2026} gives rationality. Thus \(U\) is strongly rational and \(\Rep(U)\) is modular \emph{before} module exhaustion is used.
\end{enumerate}

\section*{Exactly twelve simple modules}
The six concrete lattice constructions are
\[
\begin{array}{c|rrrrrr}
&U&B&P&Q&X_+&X_-\\\hline
h&0&1&1/3&2/3&1&1\\
\dim\mathrm{top}&1&80&3&27&81&81.
\end{array}
\]
Zero modes, positivity, and duality intertwiners establish that these are distinct simple self-contragredient modules (\textsf{C08--C09}; Section~10 and Appendix~I). The presentation allows at most twelve simples.

The conformal class has an explicit representative in \(S_K\). Its exact cocenter identity, modulo the linear span \([S_K,S_K]\) of commutators, determines by traces the six additional possible lowest weights:
\[
\frac5{13},\quad\frac6{13},\quad\frac7{13},\quad\frac8{13},
\quad\frac{11}{13},\quad\frac{15}{13}.
\]
Preserved affine idempotent ranks supply the calculation (\textsf{C16}; Section~11 and Subsection~K.6). Positive nonvacuum weights give positive categorical dimensions by character asymptotics; write \(d_Y\) for the categorical dimension of a simple module \(Y\).

The lattice branching identities force one additional actual simple. If only the six known self-contragredient simples existed, covariance would make the vacuum--\(B\) difference an \(S\)-eigenvector with eigenvalue \(\pm1\). The vacuum coordinate of its $S$-transform is \((1-d_B)/\sqrt{\sum_Y d_Y^2}\), whose absolute value is less than one because \(d_B>0\). This is a contradiction.

We now apply Coste and Gannon's Galois symmetry~\cite{CosteGannon1999} to the actual modular representation, using the square-Galois symmetry of Dong--Lin--Ng~\cite[Corollary~4.4]{DongLinNg2015}. An additional character \(T\)-eigenvalue has the form
\[
\zeta_3^{-1}\zeta_{13}^{a},\qquad a\in\{2,5,6,7,8,11\}.
\]
Squares of units modulo \(39\) fix \(\zeta_3^{-1}\) and act modulo \(13\) through \(\{1,3,4,9,10,12\}\), transitively on the six exponents (\textsf{C16--C17}; Subsection~11.5). One new eigenvalue forces six new actual simples. This is cyclotomic modular-Galois symmetry after rationality, distinct from conjugating the coefficients in \(K\). The upper bound gives exhaustion and survival of all twelve blocks. Hence
\[
S_K\otimes_K\mathbb C\cong D,\qquad
\dim D=546,\qquad \dim A(U)=36{,}625.
\]

\section*{The categorical identification}
The two isotropic lines of the discriminant form of \(E_6\oplus A_2\) give holomorphic \(E_8\) lattice extensions. Their actual vertex products supply algebra objects
\[
A_\pm=U\oplus B\oplus2X_\pm\quad\text{in}\quad\CC=\Rep(U).
\]
With the extension and separability hypotheses established, these are Lagrangian algebras. Davydov--M\"uger--Nikshych--Ostrik's center characterization~\cite[Corollary~4.1]{DMNO2013} gives
\[
\CC\simeq_{\mathrm{br}} Z(\DD_\pm),\qquad
\DD_\pm=\CC_{A_\pm}.
\]
The proposed reconstruction identifies six simple objects and the fusion rules
\[
\alpha^3=\one,\qquad \alpha\rho=\rho\alpha^2,\qquad
\rho^2=\one\oplus\rho\oplus\alpha\rho\oplus\alpha^2\rho.
\]
The finite-group and Weil-model checks, together with the modular and induction arithmetic, are recorded under \textsf{C18--C19} (Subsections~12.3--12.4 and Appendix~O). Proposition~12.5 and Appendices~M--N derive the normalized Evans--Gannon actions~\cite[Theorem~1]{EvansGannon2017} from the actual categorical splitting maps and prove tensor equivalence with the reconstructed category. The proof uses coefficient comparisons in the categorical Hom spaces and an expansion of all word morphisms to establish full faithfulness. The finite checks (\textsf{C20--C21}) verify the supplied classification equations, normalization scalars, and normalized action formulas; the passage from the actual VOA category to these formulas is proved in the manuscript. Positivity then selects \(\HH_2\) or \(\HH_3\) (Corollary~12.6 and Appendix~P; \cite[Theorem~3(ii)]{EvansGannon2017}), both in the Haagerup Morita class by Grossman--Snyder~\cite[Theorems~1.4--1.5]{GrossmanSnyder2012}.

Morita equivalence then gives the braided equivalence. The specified braiding orientation and positive spherical structures make it ribbon (Section~13, Corollary~13.2, and Appendix~Q):
\[
\boxed{\Rep(U)\simeq_{\mathrm{ribbon}}Z(\HH_3).}
\]

\bibliographystyle{plain}
\bibliography{references}

\end{document}